\documentclass[twoside,11pt]{article}

\usepackage{amsmath,amssymb, amsthm}
\usepackage{graphicx}

\newtheorem{thm}{Theorem}
\newtheorem{cor}{Corollary}
\newtheorem{lem}{Lemma}

\theoremstyle{definition}
\newtheorem{defn}{Definition}
\theoremstyle{remark}
\newtheorem{rem}{Remark}
\numberwithin{equation}{section}
\begin{document}
\pagestyle{myheadings} \markboth{ \rm \centerline {Bogdan Szal}} 
{\rm \centerline { On $L$-convergence of trigonometric series}}

\begin{titlepage}
\title{\bf {On $L$-convergence of trigonometric series}}
\author {Bogdan Szal \\
{\small University of Zielona G\'{o}ra,}\\
{\small Faculty of Mathematics, Computer Science and Econometrics,}\\
{\small 65-516 Zielona G\'{o}ra, ul. Szafrana 4a, Poland} \\ 
{\small e-mail: B.Szal @wmie.uz.zgora.pl}}
\end{titlepage}

\date{}
\maketitle
\begin{abstract} In the present paper we consider the trigonometric series with $\left( \beta
,r\right) -$general monotone and $\left( \beta ,r\right) $-rest bounded
variation coefficients. Necessary and sufficien conditions of $L$%
-convergence for such series are obtained in terms of the coefficients.
\end{abstract}



\maketitle

\section{Introduction}

Let $L\equiv L_{2\pi }$ be the space of all integrable functions $f$ of
period $2\pi $ with the norm%
\begin{equation*}
\left\Vert f\right\Vert :=\frac{1}{2\pi }\int\limits_{-\pi }^{\pi
}\left\vert f\left( x\right) \right\vert dx.
\end{equation*}%
We are interesting in the trigonometric series of the form%
\begin{equation}
\sum\limits_{k=-\infty }^{\infty }c_{k}e^{ikx}.  \label{0}
\end{equation}%
In special cases, if $c_{k}=c_{-k}=\frac{a_{k}}{2}$ then we obtain the
cosine series%
\begin{equation}
\frac{a_{0}}{2}+\sum\limits_{k=1}^{\infty }a_{k}\cos kx  \label{1}
\end{equation}%
or, when $c_{k}=-c_{-k}=-i\frac{b_{k}}{2},$ we get the sine series%
\begin{equation}
\sum\limits_{k=1}^{\infty }b_{k}\sin kx.  \label{2}
\end{equation}%
We define by $h\left( x\right) $, $f\left( x\right) $ and $g\left( x\right) $
the sums of the series (\ref{0}), (\ref{1}) and (\ref{2}), respectively,
when the series are convergent at point $x$. Moreover, denote by $%
S_{n}\left( h,x\right) $, $S_{n}\left( f,x\right) $ and $S_{n}\left(
g,x\right) $ the partial sum of these series, respectively.

We note that in the case when the sequences of the coefficients of series (%
\ref{1}) and (\ref{2}) are nonnegative monotone decreasing sequences ($%
\left( a_{n}\right) ,\left( b_{n}\right) \in M$) one has $f,g\in L$ iff
series (\ref{1}) and (\ref{2}) are the Fourier series of $f$ and $g$,
respectively. Moreover, the condition $\sum\limits_{k=1}^{\infty }\frac{%
a_{k}}{k}<\infty $ (or $\sum\limits_{k=1}^{\infty }\frac{b_{k}}{k}<\infty $%
) guarantees $f\in L$ (or $g\in L$).

Also, it is clear that if series (\ref{1}) is the Fourier series then the
condition%
\begin{equation*}
\left\Vert V_{n}\left( f\right) -S_{n}\left( f\right) \right\Vert =o\left(
1\right) ,
\end{equation*}%
where $V_{n}\left( f,x\right) =\frac{1}{n+1}\sum\limits_{k=0}^{n}S_{k}%
\left( f,x\right) $, is equivalent to convergence of $S_{n}\left( f,x\right) 
$ in $L$. The same results holds for series (\ref{2}) as well. The following
theorem provides a criterion of the convergence of $S_{n}\left( f,x\right) $
in the terms of the coefficients of (\ref{1}) and (\ref{2}) (see \cite{8})

\begin{thm}
Suppose that $\left( a_{n}\right) \in M$ and $a_{n}\rightarrow 0$. Let $f\in
L$, then%
\begin{equation}
\left\Vert f-S_{n}\left( f\right) \right\Vert =o\left( 1\right) \text{ \ \
iff \ \ }a_{n}\ln n=o\left( 1\right) .  \label{3}
\end{equation}%
The same results holds for series (\ref{2}) as well.
\end{thm}

Naturally, one would ask if the monotone condition of the coefficients can
be weakened? Indeed, many results have appeared with more general conditions
in place of the monotone condition. For example, Garrett, Rees, Stanojevi%
\'{c} \cite{23} and Teljakovski\u{\i}, Fomin \cite{5} proved that (\ref{3})
holds true for any quasimonotone sequence ($\left( a_{n}\right) ,\left(
b_{n}\right) \in QM$). Later, Stanojevic \cite{16} and Xie, Zhou \cite{7}
introduced $O$-regularly quasimonotone coefficients ($\left( a_{n}\right)
,\left( b_{n}\right) \in ORVQM$). Here,%
\begin{equation*}
QM=\left\{ \left( a_{n}\right) \in 
\mathbb{R}
_{+}:\exists \tau >0\text{ such that }\left( n^{-\tau }a_{n}\right)
\downarrow \right\}
\end{equation*}%
and%
\begin{equation*}
ORVQM=\left\{ \left( a_{n}\right) \in 
\mathbb{R}
_{+}:\exists \left( \lambda _{n}\right) \uparrow ,\lambda _{2n}\leq C\lambda
_{n}\text{ such that }\left( \lambda _{n}^{-1}a_{n}\right) \downarrow
\right\} .
\end{equation*}%
Recently, Leindler (see \cite{16}, for example) introduced the class%
\begin{equation*}
RBVS=\left\{ \left( a_{n}\right) \in 
\mathbb{C}
:\sum\limits_{n=m}^{\infty }\left\vert a_{n}-a_{n+1}\right\vert \leq
C\left\vert a_{m}\right\vert \text{ for all }m\in 
\mathbb{N}
\right\} ,
\end{equation*}%
which possess many good properties of monotone sequences, and have been used
to generalize many classical results in Fourier analysis. However Leindler
also pointed out that $RBVS$ and $QM$ cannot contain each other.

Recent investigations on $L$-convergence problem can be found in e.g. \cite%
{1, 19, 20, 21, 22, 2, 3, 4, 10, 6, 17, 18}. For example, the following
classes of coefficients were studied.

The class of general monotone coefficients $GM$, is defined as%
\begin{equation*}
GM=\left\{ \left( a_{n}\right) \in 
\mathbb{C}
:\sum\limits_{n=m}^{2m-1}\left\vert a_{n}-a_{n+1}\right\vert \leq
C\left\vert a_{m}\right\vert \text{ for all }m\in 
\mathbb{N}
\right\} .
\end{equation*}%
It turns out that for the series with $GM$-coefficients one can prove three
convergence criteria for trigonometric series in $L^{p}$: for $p=\infty $, $%
p=1$ and $1<p<\infty $ (see \cite{10}).

The $GBVS$ \cite{2} and $NBVS$ \cite{14} classes are defined as follow:

\begin{equation*}
GBVS=\left\{ \left( a_{n}\right) \in 
\mathbb{C}
\right.
\end{equation*}%
\begin{equation*}
\left. \sum\limits_{n=m}^{2m-1}\left\vert a_{n}-a_{n+1}\right\vert \leq C%
\underset{m\leq n\leq N+m}{\max }\left\vert a_{n}\right\vert \text{ for some
integer }N\text{ and all }m\in 
\mathbb{N}
\right\}
\end{equation*}%
and%
\begin{equation*}
NBVS=\left\{ a_{n}\in 
\mathbb{C}
:\sum\limits_{n=m}^{2m-1}\left\vert a_{n}-a_{n+1}\right\vert \leq C\left(
\left\vert a_{m}\right\vert +\left\vert a_{2m}\right\vert \right) \text{ for
all }m\in 
\mathbb{N}
\right\} .
\end{equation*}

For the above mentioned classes the following embedding relations are true:%
\begin{equation*}
M\varsubsetneq QM\varsubsetneq ORVQM\varsubsetneq GM\varsubsetneq GBVS\cup
NBVS
\end{equation*}%
and%
\begin{equation*}
M\varsubsetneq RBVS\varsubsetneq GM\varsubsetneq GBVS\cup NBVS.
\end{equation*}%
For the $GBVS$ and $NBVS$ classes criterion (\ref{3}) was proved in \cite{2}
and \cite{14}, respectively.

Moreover, for a more general class%
\begin{equation*}
MVBV=\left\{ a_{n}\in 
\mathbb{C}
:\sum\limits_{n=m}^{2m-1}\left\vert a_{n}-a_{n+1}\right\vert \leq
C\sum\limits_{n=\left[ m/c\right] }^{\left[ cm\right] }\frac{\left\vert
a_{n}\right\vert }{n}\text{ for some }c>1\text{ and all }m\in 
\mathbb{N}
\right\}
\end{equation*}%
criterion (\ref{3}) was considered in \cite{17}

In \cite{10, 11, 12} Tikhonov defined the class of $\beta -$general monotone
sequences as follows:

\begin{defn}
Let $\beta :=\left( \beta _{n}\right) $ be a nonnegative sequence. The
sequence of complex numbers $a:=\left( a_{n}\right) $ is said to be $\beta -$%
general monotone, or $a\in GM\left( \beta \right) $, if the relation%
\begin{equation*}
\sum\limits_{n=m}^{2m-1}\left\vert a_{n}-a_{n+1}\right\vert \leq C\beta _{m}
\end{equation*}%
holds for all $m\in 
\mathbb{N}
$.
\end{defn}

In the paper \cite{24} Tikhonov considered the following examples of the
sequences $\beta _{n}:$

(i) $_{1}\beta _{n}=\left\vert a_{n}\right\vert ,$

(ii) $_{2}\beta _{n}=\sum\limits_{k=n}^{n+N}\left\vert a_{k}\right\vert $
for some integer $N$,

(iii) $_{3}\beta _{n}=\sum\limits_{\nu =0}^{N}\left\vert a_{c^{\nu
}n}\right\vert $ for some integers $N$ and $c>1$,

(iv) $_{4}\beta _{n}=\left\vert a_{n}\right\vert +\sum\limits_{k=n+1}^{%
\left[ cn\right] }\frac{\left\vert a_{k}\right\vert }{k}$ for some $c>1$,

(v) $_{5}\beta _{n}=\sum\limits_{k=\left[ n/c\right] }^{\left[ cn\right] }%
\frac{\left\vert a_{k}\right\vert }{k}$ for some $c>1$,

(vi) $_{6}\beta _{n}=\frac{1}{\ln n}\underset{m\geq \left[ n/c\right] }{\max 
}\left( \frac{\ln m}{m}\sum\limits_{k=m}^{2m}\left\vert a_{k}\right\vert
\right) $ for some $c>1$.

We know that (see \cite{12} and \cite{24})%
\begin{equation*}
GM\left( _{1}\beta +_{2}\beta +_{3}\beta +_{4}\beta +_{5}\beta \right)
\equiv GM\left( _{5}\beta \right) \varsubsetneq GM\left( _{6}\beta \right) .
\end{equation*}

We also note that $GM\left( _{1}\beta \right) =GM$, $GM\left( _{2}\beta
\right) =GBVS$ and $NBVS\subseteq GM\left( _{3}\beta \right) $.

We write $I_{1}\ll I_{2}$ if there exists a positive constant $K$ such that $%
I_{1}\leq KI_{2}$.

In order to formulate our new results we define another classes of sequences.

\begin{defn}
\cite{13}Let $\beta :=\left( \beta _{n}\right) $ be a nonnegative sequence
and $r$ a natural number. The sequence of complex numbers $a:=\left(
a_{n}\right) $ is said to be $\left( \beta ,r\right) -$general monotone, or $%
a\in GM\left( \beta ,r\right) $, if the relation%
\begin{equation*}
\sum\limits_{n=m}^{2m-1}\left\vert a_{n}-a_{n+r}\right\vert \leq C\beta _{m}
\end{equation*}%
holds for all $m\in 
\mathbb{N}
$.
\end{defn}

\begin{defn}
Let $\beta :=\left( \beta _{n}\right) $ be a nonnegative sequence and $r$ a
natural number. The sequence of complex numbers $a:=\left( a_{n}\right) $ is
said to be $\left( \beta ,r\right) -$rest bounded variation sequence, or $%
a\in RBVS\left( \beta ,r\right) $, if the relation%
\begin{equation*}
\sum\limits_{n=m}^{\infty }\left\vert a_{n}-a_{n+r}\right\vert \leq C\beta
_{m}
\end{equation*}%
holds for all $m\in 
\mathbb{N}
$.
\end{defn}

It is clear that $RBVS\left( \beta ,r\right) \subseteq GM\left( \beta
,r\right) $ for all $r\in 
\mathbb{N}
$. Moreover, $GM\left( \beta ,1\right) \equiv GM\left( \beta \right) $ and $%
RBVS\left( \beta ,1\right) \equiv RBVS\left( \beta \right) $ ($RBVS\equiv
RBVS\left( _{1}\beta \right) )$. The next embedding relations are formulated
in the following remarks:

\begin{rem}
Let $r$ be a natural number such that $r=p\cdot q$, where $p,q\in 
\mathbb{N}
$. If a nonnegative sequence $\beta :=\left( \beta _{n}\right) $ is such that%
\begin{equation}
\sum\limits_{i=0}^{p-1}\beta _{n+i\cdot q}\ll \beta _{n}  \label{4}
\end{equation}%
for all $n$, then%
\begin{equation*}
GM\left( \beta ,q\right) \subseteq GM\left( \beta ,r\right) .
\end{equation*}
\end{rem}

\begin{rem}
Let $q,r\in 
\mathbb{N}
$ and $q\mid r$. Then%
\begin{equation*}
RBVS\left( \beta ,q\right) \subseteq RBVS\left( \beta ,r\right)
\end{equation*}
\end{rem}

It is clear that the sequences $\left( _{5}\beta _{n}\right) $ and $\left(
_{6}\beta _{n}\right) $ satisfy the condition (\ref{4}). Thus in special
case $q=1,$ from the before remark we obtain the following embedding
relations:%
\begin{equation*}
GM\left( _{5}\beta \right) \equiv GM\left( _{5}\beta ,1\right) \subseteq
GM\left( _{5}\beta ,r\right)
\end{equation*}%
and%
\begin{equation}
GM\left( _{6}\beta \right) \equiv GM\left( _{6}\beta ,1\right) \subseteq
GM\left( _{6}\beta ,r\right)  \label{5}
\end{equation}%
for all $r\in 
\mathbb{N}
$.

In this note we shall present the properties of the classes $GM\left( \beta
,r\right) $ and $RBVS\left( \beta ,r\right) .$ Moreover, we generalize and
extend to the class $GM\left( \beta ,r\right) $ and the class $RBVS\left(
\beta ,r\right) $ the Tikhonov results, which are included in \cite{24}.

\section{Main results}

The following results are true:

\begin{thm}
Let $c=\left( c_{n}\right) \in GM\left( \beta ,r\right) $, where $r\in 
\mathbb{N}
$ and a nonnegative sequence $\beta =\left( \beta _{n}\right) $ satisfies%
\begin{equation}
\sum_{k=\left[ n/2\right] }^{n}\beta _{k}\ll \sum\limits_{k=\left[ n/\gamma %
\right] }^{\left[ \gamma n\right] }\left\vert c_{k}\right\vert  \label{m1}
\end{equation}%
for some $\gamma >1$. If 
\begin{equation*}
\left\Vert V_{n}\left( h\right) -S_{n}\left( h\right) \right\Vert =o\left(
1\right)
\end{equation*}%
then $\left\vert c_{n}\right\vert \ln n=o\left( 1\right) $.
\end{thm}

\begin{thm}
Let $a=\left( a_{n}\right) \in GM\left( \beta ,2\right) $, where a
nonnegative sequence $\beta =\left( \beta _{n}\right) $ satisfies%
\begin{equation*}
\sum\limits_{k=\left[ n/2\right] }^{2n-1}\frac{\beta _{k}+\beta
_{2k}+\left\vert a_{k}\right\vert +\left\vert a_{k+1}\right\vert }{2k-n+2}%
=o\left( 1\right) .
\end{equation*}%
Then 
\begin{equation}
\left\Vert V_{n}\left( f\right) -S_{n}\left( f\right) \right\Vert =o\left(
1\right)  \label{m2a}
\end{equation}%
holds.
\end{thm}

Theorem 3 implies the following results, immediately.

\begin{cor}
Let $a=\left( a_{n}\right) \in GM\left( \beta ,2\right) $, such that%
\begin{equation*}
\left( \beta _{n}+\beta _{2n}+\left\vert a_{n}\right\vert +\left\vert
a_{n+1}\right\vert \right) \ln n=o\left( 1\right) .
\end{equation*}%
Then (\ref{m2a}) holds true.
\end{cor}

In particular, from above Corollary 1 and Theorem 2 we can derive the
following remark.

\begin{rem}
Let $a=\left( a_{n}\right) \in GM\left( _{6}\beta ,2\right) $. Then%
\begin{equation*}
\left\Vert V_{n}\left( f\right) -S_{n}\left( f\right) \right\Vert =o\left(
1\right) \text{ \ \ iff \ \ }\left\vert a_{n}\right\vert \ln n=o\left(
1\right)
\end{equation*}
\end{rem}

\begin{rem}
If we confine our attention to the class $GM\left( _{6}\beta \right) $ then
by (\ref{5}) the Tikhonov result ( see \cite[Corollary 3.3.1]{24} ) follows
from Remark 3.
\end{rem}

\begin{thm}
Let $a=\left( a_{n}\right) \in RBVS\left( \beta ,r\right) $, where $r\geq 3$%
. If a nonnegative sequence $\beta =\left( \beta _{n}\right) $ satisfies%
\begin{equation}
\sum\limits_{k=1}^{\infty }\frac{\beta _{k}}{k}<\infty  \label{m3}
\end{equation}%
then%
\begin{equation}
\left\Vert S_{n}\left( f\right) -f\right\Vert =o\left( 1\right) .  \label{m4}
\end{equation}
\end{thm}

\begin{thm}
Let $b=\left( b_{n}\right) \in RBVS\left( \beta ,r\right) $, where $r\geq 2$
and a nonnegative sequence $\beta =\left( \beta _{n}\right) $ satisfies (\ref%
{m3}). Then%
\begin{equation}
\left\Vert S_{n}\left( g\right) -g\right\Vert =o\left( 1\right) .  \label{m5}
\end{equation}
\end{thm}

\begin{rem}
(i) There exists a sequence $\left( a_{n}\right) \in RBVS\left( _{5}\beta
,3\right) $ with the property $a_{n}\ln n=o\left( 1\right) $ such that the
series (\ref{1}) is not convergent in $L$-norm.

(ii) There exists a sequence $\left( b_{n}\right) \in RBVS\left( _{5}\beta
,2\right) $ with the property $b_{n}\ln n=o\left( 1\right) $ such that the
series (\ref{2}) is not convergent in $L$-norm.
\end{rem}

\begin{rem}
For any $r\geq 3$ there exists a sequence $d=\left( d_{n}\right) \in
RBV\left( _{6}\beta ,r\right) $ with the property $\sum\limits_{k=1}^{%
\infty }\frac{\beta _{k}}{k}<\infty $, which does not belong to the class $%
GM\left( _{6}\beta ,2\right) .$
\end{rem}

Combining Theorem 3, 4 and 5 we obtain the following assertion.

\begin{rem}
Let $c=\left( c_{n}\right) \in RBV\left( \beta ,r\right) $, where $r\geq 2$.
If a nonnegative sequence $\beta =\left( \beta _{n}\right) $ satisfies%
\begin{equation*}
\sum\limits_{k=1}^{\infty }\frac{\beta _{k}}{k}<\infty
\end{equation*}%
then%
\begin{equation*}
\left\Vert S_{n}\left( h\right) -h\right\Vert =o\left( 1\right) .
\end{equation*}
\end{rem}

\section{Lemmas}

\begin{lem}
Let $r\in 
\mathbb{N}
$, $l\in 
\mathbb{Z}
$ and $a:=\left( a_{n}\right) \in 
\mathbb{C}
$. If $x\neq \frac{2l\pi }{r}$, then for all $m\geq n$%
\begin{equation*}
\sum\limits_{k=n}^{m}a_{k}\cos kx=\frac{1}{2\sin \left( rx/2\right) }%
\left\{ \sum\limits_{k=n}^{m}\left( a_{k}-a_{k+r}\right) \sin \left( k+%
\frac{r}{2}\right) x\right.
\end{equation*}%
\begin{equation}
\left. +\sum\limits_{k=m+1}^{m+r}a_{k}\sin \left( k-\frac{r}{2}\right)
x-\sum\limits_{k=n}^{n+r-1}a_{k}\sin \left( k-\frac{r}{2}\right) x\right\}
\label{l01}
\end{equation}%
and%
\begin{equation*}
\sum\limits_{k=n}^{m}a_{k}\sin kx=\frac{-1}{2\sin \left( rx/2\right) }%
\left\{ \sum\limits_{k=n}^{m}\left( a_{k}-a_{k+r}\right) \cos \left( k+%
\frac{r}{2}\right) x\right.
\end{equation*}%
\begin{equation}
\left. +\sum\limits_{k=m+1}^{m+r}a_{k}\cos \left( k-\frac{r}{2}\right)
x-\sum\limits_{k=n}^{n+r-1}a_{k}\cos \left( k-\frac{r}{2}\right) x\right\} .
\label{l02}
\end{equation}
\end{lem}

\begin{proof}
We start with proof of the equality (\ref{l01}). An elementary calculation
gives%
\begin{equation*}
\sum\limits_{k=n}^{m}a_{k}\cos kx=\sum\limits_{k=n}^{m}\left(
a_{k}-a_{k+r}\right) \cos kx+\sum\limits_{k=n}^{m}a_{k+r}\cos kx
\end{equation*}%
\begin{equation}
=\sum\limits_{k=n}^{m}\left( a_{k}-a_{k+r}\right) \cos kx+\cos
rx\sum\limits_{k=n+r}^{m+r}a_{k}\cos kx+\sin
rx\sum\limits_{k=n+r}^{m+r}a_{k}\sin kx.  \label{l1}
\end{equation}%
On the other hand%
\begin{equation*}
\sum\limits_{k=n}^{m}a_{k}\sin kx=\sum\limits_{k=n}^{m}\left(
a_{k}-a_{k+r}\right) \sin kx+\sum\limits_{k=n}^{m}a_{k+r}\sin kx
\end{equation*}%
\begin{equation}
=\sum\limits_{k=n}^{m}\left( a_{k}-a_{k+r}\right) \sin kx+\cos
rx\sum\limits_{k=n+r}^{2n+r-1}a_{k}\sin kx-\sin
rx\sum\limits_{k=n+r}^{2n+r-1}a_{k}\cos kx.  \label{l2}
\end{equation}%
Hence%
\begin{equation*}
\left( 1-\cos rx\right) \sum\limits_{k=n+r}^{m+r}a_{k}\sin
kx=\sum\limits_{k=n}^{m}\left( a_{k}-a_{k+r}\right) \sin kx-\sin
rx\sum\limits_{k=n+r}^{m+r}a_{k}\cos kx
\end{equation*}%
\begin{equation*}
-\sum\limits_{k=n}^{n+r-1}a_{k}\sin kx+\sum\limits_{k=m+1}^{m+r}a_{k}\sin
kx.
\end{equation*}%
Therefore, if $x\neq \frac{2l\pi }{r}$, then%
\begin{equation*}
\sum\limits_{k=n+r}^{m+r}a_{k}\sin kx=\frac{1}{2\sin ^{2}\left( rx/2\right) 
}\left\{ \sum\limits_{k=n}^{m}\left( a_{k}-a_{k+r}\right) \sin kx\right.
\end{equation*}%
\begin{equation*}
\left. -\sin rx\sum\limits_{k=n+r}^{m+r}a_{k}\cos
kx-\sum\limits_{k=n}^{n+r-1}a_{k}\sin
kx+\sum\limits_{k=m+1}^{m+r}a_{k}\sin kx\right\} .
\end{equation*}%
Putting this to (\ref{l1}) we get%
\begin{equation*}
\sum\limits_{k=n}^{m}a_{k}\cos kx=\sum\limits_{k=n}^{m}\left(
a_{k}-a_{k+r}\right) \cos kx+\cos rx\sum\limits_{k=n+r}^{m+r}a_{k}\cos kx
\end{equation*}%
\begin{equation*}
+\frac{\cos \left( rx/2\right) }{\sin \left( rx/2\right) }%
\sum\limits_{k=n}^{m}\left( a_{k}-a_{k+r}\right) \sin kx-2\cos ^{2}\frac{rx%
}{2}\sum\limits_{k=n+r}^{m+r}a_{k}\cos kx
\end{equation*}%
\begin{equation*}
-\frac{\cos \left( rx/2\right) }{\sin \left( rx/2\right) }%
\sum\limits_{k=n}^{n+r-1}a_{k}\sin kx+\frac{\cos \left( rx/2\right) }{\sin
\left( rx/2\right) }\sum\limits_{k=m+1}^{m+r}a_{k}\sin kx
\end{equation*}%
\begin{equation*}
=\frac{1}{\sin \left( rx/2\right) }\sum\limits_{k=n}^{2n-1}\left(
a_{k}-a_{k+r}\right) \sin \left( k+\frac{r}{2}\right)
-\sum\limits_{k=n+r}^{m+r}a_{k}\cos kx
\end{equation*}%
\begin{equation*}
-\frac{\cos \left( rx/2\right) }{\sin \left( rx/2\right) }%
\sum\limits_{k=n}^{n+r-1}a_{k}\sin kx+\frac{\cos \left( rx/2\right) }{\sin
\left( rx/2\right) }\sum\limits_{k=m+1}^{m+r}a_{k}\sin kx.
\end{equation*}%
Thus%
\begin{equation*}
2\sum\limits_{k=n}^{m}a_{k}\cos kx=\frac{1}{\sin \left( rx/2\right) }%
\sum\limits_{k=n}^{m}\left( a_{k}-a_{k+r}\right) \sin \left( k+\frac{r}{2}%
\right)
\end{equation*}%
\begin{equation*}
-\frac{\cos \left( rx/2\right) }{\sin \left( rx/2\right) }%
\sum\limits_{k=n}^{n+r-1}a_{k}\sin kx+\frac{\cos \left( rx/2\right) }{\sin
\left( rx/2\right) }\sum\limits_{k=m+1}^{m+r}a_{k}\sin kx
\end{equation*}%
\begin{equation*}
+\sum\limits_{k=n}^{n+r-1}a_{k}\cos kx-\sum\limits_{k=m+1}^{m+r}a_{k}\cos
kx
\end{equation*}%
\begin{equation*}
=\frac{1}{\sin \left( rx/2\right) }\left\{ \sum\limits_{k=n}^{m}\left(
a_{k}-a_{k+r}\right) \sin \left( k+\frac{r}{2}\right) x\right.
\end{equation*}%
\begin{equation*}
\left. +\sum\limits_{k=m+1}^{m+r}a_{k}\sin \left( k-\frac{r}{2}\right)
x-\sum\limits_{k=n}^{n+r-1}a_{k}\sin \left( k-\frac{r}{2}\right) x\right\}
\end{equation*}%
and (\ref{l01}) holds.

Now we show that (\ref{l02}) is true, too. By (\ref{l1}) we get%
\begin{equation*}
\left( 1-\cos rx\right) \sum\limits_{k=n+r}^{m+r}a_{k}\cos
kx=\sum\limits_{k=n}^{m}\left( a_{k}-a_{k+r}\right) \cos kx+\sin
rx\sum\limits_{k=n+r}^{m+r}a_{k}\sin kx
\end{equation*}%
\begin{equation*}
-\sum\limits_{k=n}^{n+r-1}a_{k}\cos kx+\sum\limits_{k=m+1}^{m+r}a_{k}\cos
kx.
\end{equation*}%
Therefore, if $x\neq \frac{2l\pi }{r}$, then%
\begin{equation*}
\sum\limits_{k=n+r}^{m+r}a_{k}\cos kx=\frac{1}{2\sin ^{2}\left( rx/2\right) 
}\left\{ \sum\limits_{k=n}^{m}\left( a_{k}-a_{k+r}\right) \cos kx\right.
\end{equation*}%
\begin{equation*}
\left. +\sin rx\sum\limits_{k=n+r}^{m+r}a_{k}\sin
kx-\sum\limits_{k=n}^{n+r-1}a_{k}\cos
kx+\sum\limits_{k=2n}^{2n+r-1}a_{k}\cos kx\right\} .
\end{equation*}%
Putting this to (\ref{l2}) and making the same calculation as before we
obtain that (\ref{l02}) holds, too.

The proof is complete.
\end{proof}

\section{Proofs of the main results}

In this section we shall prove our theorems and remarks.

\subsection{Proof of Remark 1}

Let $r=p\cdot q$ and $\left( a_{n}\right) \in GM\left( \beta ,q\right) $.
Then using (\ref{4}) we get%
\begin{equation*}
\sum\limits_{k=n}^{2n-1}\left\vert a_{k}-a_{k+r}\right\vert
=\sum\limits_{k=n}^{2n-1}\left\vert \sum\limits_{l=0}^{p-1}\left(
a_{k+l\cdot q}-a_{k+\left( l+1\right) \cdot q}\right) \right\vert
\end{equation*}%
\begin{equation*}
\leq \sum\limits_{k=n}^{2n-1}\sum\limits_{l=0}^{p-1}\left\vert a_{k+l\cdot
q}-a_{k+\left( l+1\right) \cdot q}\right\vert
=\sum\limits_{l=0}^{p-1}\sum\limits_{k=n+l\cdot q}^{2n+l\cdot
q-1}\left\vert a_{k}-a_{k+q}\right\vert \ll \sum\limits_{l=0}^{p-1}\beta
_{n+l\cdot q}\ll \beta _{n}
\end{equation*}%
whence $\left( a_{n}\right) GM\left( \beta ,r\right) $. $\square $

\subsection{Proof of Remark 1}

Let $q,r\in 
\mathbb{N}
$ and $q\mid r$. Then there exists a natural number $p$ such that $r=p\cdot
q.$ Supposing $\left( a_{n}\right) \in RBVS\left( \beta ,q\right) $ we have
for all $n$%
\begin{equation*}
\sum\limits_{k=n}^{\infty }\left\vert a_{k}-a_{k+r}\right\vert
=\sum\limits_{k=n}^{\infty }\left\vert \sum\limits_{l=0}^{p-1}\left(
a_{k+l\cdot q}-a_{k+\left( l+1\right) \cdot q}\right) \right\vert
\end{equation*}%
\begin{equation*}
\leq \sum\limits_{k=n}^{\infty }\sum\limits_{l=0}^{p-1}\left\vert
a_{k+l\cdot q}-a_{k+\left( l+1\right) \cdot q}\right\vert
=\sum\limits_{l=0}^{p-1}\sum\limits_{k=n+l\cdot q}^{\infty }\left\vert
a_{k}-a_{k+q}\right\vert \leq p\sum\limits_{k=n}^{\infty }\left\vert
a_{k}-a_{k+q}\right\vert \ll \beta _{n}.
\end{equation*}%
Hence $\left( a_{n}\right) \in RBVS\left( \beta ,r\right) $ and this ends
our proof. $\square $

\subsection{Proof of Theorem 2}

Let $\left( c_{n}\right) \in GM\left( \beta ,r\right) $, where $r\geq 1$. We
will show that 
\begin{equation*}
n\left\vert c_{n}\right\vert \ll \sum\limits_{k=\left[ n/\gamma \right] }^{%
\left[ \gamma n\right] }\left\vert c_{k}\right\vert \text{, \ \ }\gamma >1.
\end{equation*}%
If $n\leq r,$ then the inequality obviously holds.

Now, let $n>r$. For $j=n+1,n+2,...,2n$ we get%
\begin{equation*}
\sum\limits_{k=n}^{j-1}\left\vert c_{k}-c_{k+r}\right\vert \geq
\sum\limits_{k=n}^{j-1}\left\vert \left\vert c_{k}\right\vert -\left\vert
c_{k+r}\right\vert \right\vert \geq \left\vert
\sum\limits_{k=n}^{j-1}\left( \left\vert c_{k}\right\vert -\left\vert
c_{k+r}\right\vert \right) \right\vert \geq
\sum\limits_{k=n}^{j-1}\left\vert c_{k}\right\vert
-\sum\limits_{k=n+r}^{j+r-1}\left\vert c_{k}\right\vert
\end{equation*}%
and for $j\geq n+r+1$ we obtain%
\begin{equation*}
\left\vert c_{n}\right\vert \leq \sum\limits_{k=n}^{n+r-1}\left\vert
c_{k}\right\vert \leq \sum\limits_{k=n}^{j-1}\left\vert
c_{k}-c_{k+r}\right\vert +\sum\limits_{k=j}^{j+r-1}\left\vert
c_{k}\right\vert
\end{equation*}%
\begin{equation*}
\leq \sum\limits_{k=\left[ j/2\right] }^{2\left[ j/2\right] -1}\left\vert
c_{k}-c_{k+r}\right\vert +\sum\limits_{k=j}^{j+r-1}\left\vert
c_{k}\right\vert \ll \beta _{\left[ j/2\right] }+\sum\limits_{k=j}^{j+r-1}%
\left\vert c_{k}\right\vert .
\end{equation*}%
Summing up on $j$ and using (\ref{m1}) we get%
\begin{equation*}
n\left\vert c_{n}\right\vert =\sum\limits_{j=n+1}^{2n}\left\vert
c_{n}\right\vert =\sum\limits_{j=n+1}^{n+r}\left\vert c_{n}\right\vert
+\sum\limits_{j=n+r+1}^{2n}\left\vert c_{n}\right\vert
\end{equation*}%
\begin{equation*}
\ll r\left\vert c_{n}\right\vert +\sum\limits_{j=n+r+1}^{2n}\left( \beta _{ 
\left[ j/2\right] }+\sum\limits_{k=j}^{j+r-1}\left\vert c_{k}\right\vert
\right)
\end{equation*}%
\begin{equation*}
\leq r\left\vert c_{n}\right\vert +\sum\limits_{j=n+r+1}^{2n}\beta _{\left[
j/2\right] }+\sum\limits_{j=0}^{r-1}\sum\limits_{k=n+1+j}^{2n-1+j}\left%
\vert c_{k}\right\vert \leq \sum\limits_{j=\left[ n/2\right] }^{n}\beta
_{j}+r\sum\limits_{k=n}^{2n+r-2}\left\vert c_{k}\right\vert \ll
\sum\limits_{k=\left[ n/\gamma \right] }^{\left[ \gamma n\right]
}\left\vert c_{k}\right\vert .
\end{equation*}%
Further, applying Theorem (8.11) \cite[Chapter 7]{15} we obtain%
\begin{equation*}
\left\vert c_{n}\right\vert \ln n=\frac{\ln n}{n}\sum\limits_{k=\left[
n/\gamma \right] }^{\left[ \gamma n\right] }\left\vert c_{k}\right\vert
=o\left( 1\right) .
\end{equation*}%
The proof is complete. $\square $

\subsection{Proof of Theorem 3}

First, using inequality (see \cite{19})%
\begin{equation*}
\left\Vert V_{n}\left( f\right) -S_{n}\left( f\right) \right\Vert \ll \frac{1%
}{n+1}\sum\limits_{j=1}^{n}\left\Vert S_{j}\left( f\right) -S_{\left[ j/2%
\right] }\right\Vert +\underset{k=\left[ n/2\right] ,...,n}{\max }\left\Vert
S_{k}\left( f\right) -S_{\left[ n/2\right] }\left( f\right) \right\Vert
\end{equation*}%
we note that%
\begin{equation}
\underset{n\leq m\leq 2n-1}{\max }\left\Vert S_{m}\left( f\right)
-S_{n-1}\left( f\right) \right\Vert _{1}=o\left( 1\right)  \label{p1}
\end{equation}%
implies (\ref{m2a}). Let us now show that our condition on $\left(
a_{n}\right) $ quarantines the accuracy of (\ref{p1}). Indeed, by Lemma 1,

\begin{equation*}
\left\Vert S_{m}\left( f\right) -S_{n-1}\left( f\right) \right\Vert
=\left\Vert \sum\limits_{k=n}^{m}a_{k}\cos \left( k\cdot \right) \right\Vert
\end{equation*}%
\begin{equation*}
=\frac{1}{\pi }\int\limits_{0}^{\pi }\left\vert \frac{1}{2\sin x}\left\{
\sum\limits_{k=n}^{m}\left( a_{k}-a_{k+2}\right) \sin \left( k+1\right)
x\right. \right.
\end{equation*}%
\begin{equation*}
\left. \left. +\sum\limits_{k=m+1}^{m+2}a_{k}\sin \left( k-1\right)
x-\sum\limits_{k=n}^{n+1}a_{k}\sin \left( k-1\right) x\right\} \right\vert
dx
\end{equation*}%
\begin{equation*}
=\frac{1}{2\pi }\left\{ \sum\limits_{k=n}^{m}\left\vert
a_{k}-a_{k+2}\right\vert \int\limits_{0}^{\pi }\frac{\left\vert \sin \left(
k+1\right) x\right\vert }{\left\vert \sin x\right\vert }dx\right.
\end{equation*}%
\begin{equation*}
\left. +\sum\limits_{k=m+1}^{m+2}\left\vert a_{k}\right\vert
\int\limits_{0}^{\pi }\frac{\left\vert \sin \left( k-1\right) x\right\vert 
}{\left\vert \sin x\right\vert }dx+\sum\limits_{k=n}^{n+1}\left\vert
a_{k}\right\vert \int\limits_{0}^{\pi }\frac{\left\vert \sin \left(
k-1\right) x\right\vert }{\left\vert \sin x\right\vert }dx\right\} .
\end{equation*}%
Further, we estimate the above integrals. For $k\geq n$ we get%
\begin{equation*}
\int\limits_{0}^{\pi }\frac{\left\vert \sin \left( k+1\right) x\right\vert 
}{\left\vert \sin x\right\vert }dx
\end{equation*}%
\begin{equation*}
=\left( \int\limits_{0}^{\pi /2\left( 2k-n+1\right) }+\int\limits_{2\pi
/\left( 2k-n+1\right) }^{\pi /2}+\int\limits_{\pi /2}^{\pi -\pi /2\left(
2k-n+1\right) }+\int\limits_{\pi -\pi /2\left( 2k-n+1\right) }^{\pi
}\right) \frac{\left\vert \sin \left( k+1\right) x\right\vert }{\sin x}dx
\end{equation*}%
\begin{equation*}
=I_{1}+I_{2}+I_{3}+I_{4}.
\end{equation*}%
Using the inequalities : 
\begin{equation*}
\left\vert \sin \left( k+1\right) x\right\vert \leq \left( k+1\right)
\left\vert \sin x\right\vert ,\text{ \ \ }\sin x\geq \frac{2}{\pi }x\text{ \
\ for \ \ }x\in \left[ 0,\frac{\pi }{2}\right]
\end{equation*}%
and%
\begin{equation*}
\sin x\geq 2\left( 1-\frac{1}{\pi }x\right) \text{ \ \ for \ \ }x\in \left[ 
\frac{\pi }{2},\pi \right]
\end{equation*}%
we obtain that%
\begin{equation*}
I_{1}\leq \left( k+1\right) \int\limits_{0}^{\pi /2\left( 2k-n+1\right) }dx=%
\frac{\pi }{2}\frac{k+1}{2k-n+1}\leq \frac{\pi }{2},
\end{equation*}%
\begin{equation*}
I_{2}\leq \frac{\pi }{2}\int\limits_{2\pi /\left( 2k-n+1\right) }^{\pi /2}%
\frac{1}{x}dx=\frac{\pi }{2}\ln \left( 2k-n+1\right) ,
\end{equation*}%
\begin{equation*}
I_{3}\leq \frac{1}{2}\int\limits_{\pi /2}^{\pi -\pi /2\left( 2k-n+1\right) }%
\frac{1}{1-\frac{1}{\pi }x}dx=\frac{\pi }{2}\ln \left( 2k-n+1\right)
\end{equation*}%
and%
\begin{equation*}
I_{4}\leq \left( k+1\right) \int\limits_{\pi -\pi /2\left( 2k-n+1\right)
}^{\pi }dx=\frac{\pi }{2}\frac{k+1}{2k-n+1}\leq \frac{\pi }{2}.
\end{equation*}%
Therefore, for $k\geq n,$%
\begin{equation*}
\int\limits_{0}^{\pi }\frac{\left\vert \sin \left( k+1\right) x\right\vert 
}{\left\vert \sin x\right\vert }dx\leq \pi \left( 1+\ln \left( 2k-n+1\right)
\right) \leq 2\pi \ln \left( 2k-n+1\right) .
\end{equation*}%
Similarly we can show that%
\begin{equation*}
\int\limits_{0}^{\pi }\frac{\left\vert \sin \left( k-1\right) x\right\vert 
}{\left\vert \sin x\right\vert }dx\leq 2\pi \ln \left( 2k-n+1\right) .
\end{equation*}%
Hence%
\begin{equation*}
\left\Vert S_{m}\left( f\right) -S_{n-1}\left( f\right) \right\Vert \leq
\sum\limits_{k=n}^{m}\left\vert a_{k}-a_{k+2}\right\vert \ln \left(
2k-n+1\right)
\end{equation*}%
\begin{equation*}
+\sum\limits_{k=m+1}^{m+2}\left\vert a_{k}\right\vert \ln \left(
2k-n+1\right) +\sum\limits_{k=n}^{n+1}\left\vert a_{k}\right\vert \ln
\left( 2k-n+1\right) .
\end{equation*}%
It is clear that for $m\geq n$%
\begin{equation*}
\sum\limits_{k=n}^{m}\left\vert a_{k}-a_{k+2}\right\vert \geq \left\vert
\sum\limits_{k=n}^{n+1}\left\vert a_{k}\right\vert
-\sum\limits_{k=m+1}^{m+2}\left\vert a_{k}\right\vert \right\vert .
\end{equation*}%
Therefore, if $\left( a_{n}\right) \in GM\left( \beta ,2\right) $, then 
\begin{equation*}
\left\Vert S_{m}\left( f\right) -S_{n-1}\left( f\right) \right\Vert \leq
\sum\limits_{k=n}^{m}\left\vert a_{k}-a_{k+2}\right\vert \ln \left(
2k-n+1\right)
\end{equation*}%
\begin{equation*}
+\sum\limits_{k=m+1}^{m+2}\left\vert a_{k}\right\vert \ln \left(
2k-n+1\right) +\ln \left( n+3\right) \left( \sum\limits_{k=n}^{m}\left\vert
a_{k}-a_{k+2}\right\vert +\sum\limits_{k=m+1}^{m+2}\left\vert
a_{k}\right\vert \right)
\end{equation*}%
\begin{equation*}
\leq 2\left( \sum\limits_{k=n}^{m}\left\vert a_{k}-a_{k+2}\right\vert \ln
\left( 2k-n+3\right) +\ln \left( 2m-n+3\right)
\sum\limits_{k=m+1}^{m+2}\left\vert a_{k}\right\vert \right)
\end{equation*}%
\begin{equation*}
\ll \sum\limits_{k=n}^{m}\left\vert a_{k}-a_{k+2}\right\vert
\sum\limits_{l=\left[ n/2\right] }^{k}\frac{1}{2l-n+2}+\sum\limits_{l=%
\left[ n/2\right] }^{m}\frac{1}{2l-n+2}\sum\limits_{k=m+1}^{m+2}\left\vert
a_{k}\right\vert
\end{equation*}%
\begin{equation*}
\leq \sum\limits_{k=\left[ n/2\right] }^{m}\left\vert
a_{k}-a_{k+2}\right\vert \sum\limits_{l=\left[ n/2\right] }^{k}\frac{1}{%
2l-n+2}
\end{equation*}%
\begin{equation*}
+\sum\limits_{l=\left[ n/2\right] }^{m}\frac{1}{2l-n+2}\left(
\sum\limits_{k=l}^{l+1}\left\vert a_{k}\right\vert
+\sum\limits_{k=l}^{m}\left\vert a_{k}-a_{k+2}\right\vert \right)
\end{equation*}%
\begin{equation*}
\leq 2\sum\limits_{l=\left[ n/2\right] }^{m}\frac{1}{2l-n+2}\left(
\sum\limits_{k=l}^{l+1}\left\vert a_{k}\right\vert
+\sum\limits_{k=l}^{m}\left\vert a_{k}-a_{k+2}\right\vert \right)
\end{equation*}%
\begin{equation*}
\leq 2\sum\limits_{l=\left[ n/2\right] }^{2n-1}\frac{1}{2l-n+2}\left(
\sum\limits_{k=l}^{l+1}\left\vert a_{k}\right\vert
+\sum\limits_{k=l}^{2l-1}\left\vert a_{k}-a_{k+2}\right\vert
+\sum\limits_{k=2l}^{4l-1}\left\vert a_{k}-a_{k+2}\right\vert \right)
\end{equation*}%
\begin{equation*}
\ll \sum\limits_{k=\left[ n/2\right] }^{2n-1}\frac{\beta _{k}+\beta
_{2k}+\left\vert a_{k}\right\vert +\left\vert a_{k+1}\right\vert }{2k-n+2}.
\end{equation*}%
Thus, we obtain (\ref{p1}), which yields (\ref{m2a}) and the proof is
complete. $\square $

\subsection{Proof of Theorem 4}

If $\left( a_{n}\right) \in RBV\left( \beta ,r\right) $ then%
\begin{equation*}
\sum\limits_{k=1}^{\infty }\frac{\left\vert a_{k}\right\vert }{k}\ll
\sum\limits_{k=1}^{\infty }\frac{\beta _{k}}{k}<\infty .
\end{equation*}%
Thus $a_{n}=o\left( 1\right) $ and by (\ref{l01}) we obtain that, for $x\neq 
\frac{2l\pi }{r}$, $l\in 
\mathbb{Z}
,$%
\begin{equation*}
\sum\limits_{k=n}^{\infty }a_{k}\cos kx=\frac{1}{2\sin \left( rx/2\right) }%
\left\{ \sum\limits_{k=n}^{\infty }\left( a_{k}-a_{k+r}\right) \sin \left(
k+\frac{r}{2}\right) x\right.
\end{equation*}%
\begin{equation*}
\left. -\sum\limits_{k=n}^{n+r-1}a_{k}\sin \left( k-\frac{r}{2}\right)
x\right\} .
\end{equation*}%
Therefore%
\begin{equation*}
\left\vert f\left( x\right) -S_{n}\left( f;x\right) \right\vert =\left\vert
\sum\limits_{k=n+1}^{\infty }a_{k}\cos kx\right\vert
\end{equation*}%
\begin{equation}
\leq \frac{1}{2\left\vert \sin \left( rx/2\right) \right\vert }\left\{
\sum\limits_{k=n+1}^{\infty }\left\vert a_{k}-a_{k+r}\right\vert
+\sum\limits_{k=n+1}^{n+r}\left\vert a_{k}\right\vert \right\} \ll \frac{%
\beta _{n+1}}{\left\vert \sin \left( rx/2\right) \right\vert }.  \label{p3p1}
\end{equation}%
It is clear that for an odd $r$%
\begin{equation*}
\left\Vert S_{n}\left( f\right) -f\right\Vert =\frac{1}{\pi }%
\int\limits_{0}^{\pi }\left\vert f\left( x\right) -S_{n}\left( f;x\right)
\right\vert dx
\end{equation*}%
\begin{equation*}
=\frac{1}{\pi }\left\{ \sum\limits_{k=0}^{\left[ r/2\right]
-1}\int\limits_{\frac{2k\pi }{r}}^{\frac{2\left( k+1\right) \pi }{r}%
}\left\vert f\left( x\right) -S_{n}\left( f;x\right) \right\vert
dx+\int\limits_{\frac{2\left[ r/2\right] \pi }{r}}^{\pi }\left\vert f\left(
x\right) -S_{n}\left( f;x\right) \right\vert dx\right\}
\end{equation*}%
\begin{equation}
=\frac{1}{\pi }\left\{ \sum\limits_{k=0}^{\left[ r/2\right] }\int\limits_{%
\frac{2k\pi }{r}}^{\frac{2k\pi }{r}+\frac{\pi }{r}}\left\vert f\left(
x\right) -S_{n}\left( f;x\right) \right\vert dx+\sum\limits_{k=0}^{\left[
r/2\right] -1}\int\limits_{\frac{2k\pi }{r}+\frac{\pi }{r}}^{\frac{2\left(
k+1\right) \pi }{r}}\left\vert f\left( x\right) -S_{n}\left( f;x\right)
\right\vert dx\right\}  \label{s1}
\end{equation}%
and for an even $r$%
\begin{equation}
\left\Vert S_{n}\left( f\right) -f\right\Vert =\frac{1}{\pi }\left\{
\sum\limits_{k=0}^{\left[ r/2\right] }\left( \int\limits_{\frac{2k\pi }{r}%
}^{\frac{2k\pi }{r}+\frac{\pi }{r}}+\int\limits_{\frac{2k\pi }{r}+\frac{\pi 
}{r}}^{\frac{2\left( k+1\right) \pi }{r}}\right) \left\vert f\left( x\right)
-S_{n}\left( f;x\right) \right\vert dx\right\} .  \label{s2}
\end{equation}%
Let 
\begin{equation*}
\frac{2k\pi }{r}+\frac{\pi }{M+1}<x\leq \frac{2k\pi }{r}+\frac{\pi }{r},
\end{equation*}%
where $M:=M\left( x\right) \geq r$ and $k=0,1,...,\left[ r/2\right] -1$ if $%
r $ is an even number, and $k=0,1,...,\left[ r/2\right] $ if $r$ is an odd
number.

Then, for $n\geq M$, by (\ref{p3p1}) we get%
\begin{equation*}
\int\limits_{\frac{2k\pi }{r}+\frac{\pi }{n+1}}^{\frac{2k\pi }{r}+\frac{\pi 
}{r}}\left\vert f\left( x\right) -S_{n}\left( f;x\right) \right\vert
dx=\sum\limits_{M=r}^{n}\int\limits_{\frac{2k\pi }{r}+\frac{\pi }{M+1}}^{%
\frac{2k\pi }{r}+\frac{\pi }{M}}\left\vert f\left( x\right) -S_{n}\left(
f;x\right) \right\vert dx
\end{equation*}%
\begin{equation*}
\ll \sum\limits_{M=r}^{n}\int\limits_{\frac{2k\pi }{r}+\frac{\pi }{M+1}}^{%
\frac{2k\pi }{r}+\frac{\pi }{M}}\frac{\beta _{n+1}}{\left\vert \sin \left(
rx/2\right) \right\vert }dx.
\end{equation*}%
Using the inequality%
\begin{equation}
\frac{r}{\pi }x-2k\leq \left\vert \sin \left( \frac{r}{2}x\right)
\right\vert \text{ \ \ for \ }x\in \left[ \frac{2k\pi }{r},\frac{2k\pi }{r}+%
\frac{\pi }{r}\right]  \label{p3p2}
\end{equation}%
we have%
\begin{equation*}
\int\limits_{\frac{2k\pi }{r}+\frac{\pi }{n+1}}^{\frac{2k\pi }{r}+\frac{\pi 
}{r}}\left\vert f\left( x\right) -S_{n}\left( f;x\right) \right\vert dx\ll
\beta _{n+1}\sum\limits_{M=r}^{n}\int\limits_{\frac{2k\pi }{r}+\frac{\pi }{%
M+1}}^{\frac{2k\pi }{r}+\frac{\pi }{M}}\frac{1}{\frac{r}{\pi }x-2k}dx
\end{equation*}%
\begin{equation}
=\frac{\pi }{r}\beta _{n+1}\sum\limits_{M=r}^{n}\frac{1}{M}\ll \beta
_{n+1}\ln \left( n+1\right) .  \label{s3}
\end{equation}%
If $M\geq n+1$ then by (\ref{p3p1})%
\begin{equation*}
\left\vert \sum\limits_{k=n+1}^{\infty }a_{k}\cos kx\right\vert \leq
\left\vert \sum\limits_{k=n+1}^{M}a_{k}\cos kx\right\vert +\left\vert
\sum\limits_{k=M+1}^{\infty }a_{k}\cos kx\right\vert
\end{equation*}%
\begin{equation}
\ll \sum\limits_{k=n+1}^{M}\left\vert a_{k}\right\vert +\frac{\beta _{M+1}}{%
\left\vert \sin \left( rx/2\right) \right\vert }.  \label{p3p3}
\end{equation}%
Hence, using (\ref{p3p2}),%
\begin{equation*}
\int\limits_{\frac{2k\pi }{r}}^{\frac{2k\pi }{r}+\frac{\pi }{n+1}%
}\left\vert f\left( x\right) -S_{n}\left( f;x\right) \right\vert
dx=\sum\limits_{M=n+1}^{\infty }\int\limits_{\frac{2k\pi }{r}+\frac{\pi }{%
M+1}}^{\frac{2k\pi }{r}+\frac{\pi }{M}}\left\vert f\left( x\right)
-S_{n}\left( f;x\right) \right\vert dx
\end{equation*}%
\begin{equation*}
\ll \sum\limits_{M=n+1}^{\infty }\int\limits_{\frac{2k\pi }{r}+\frac{\pi }{%
M+1}}^{\frac{2k\pi }{r}+\frac{\pi }{M}}\left(
\sum\limits_{k=n+1}^{M}\left\vert a_{k}\right\vert +\frac{\beta _{M+1}}{%
\left\vert \sin \left( rx/2\right) \right\vert }\right) dx
\end{equation*}%
\begin{equation*}
\leq \sum\limits_{M=n+1}^{\infty }\left\{ \frac{\pi }{M\left( M+1\right) }%
\sum\limits_{k=n+1}^{M}\left\vert a_{k}\right\vert +\frac{\pi }{r}\beta
_{M+1}\int\limits_{\frac{2k\pi }{r}+\frac{\pi }{M+1}}^{\frac{2k\pi }{r}+%
\frac{\pi }{M}}\frac{1}{\frac{r}{\pi }x-2k}dx\right\}
\end{equation*}%
\begin{equation*}
\ll \sum\limits_{M=n+1}^{\infty }\frac{1}{M^{2}}\sum\limits_{k=n+1}^{M}%
\left\vert a_{k}\right\vert +\sum\limits_{M=n+1}^{\infty }\frac{\beta _{M+1}%
}{M}\leq \sum\limits_{k=n+1}^{\infty }\left\vert a_{k}\right\vert
\sum\limits_{M=k}^{\infty }\frac{1}{M^{2}}+\sum\limits_{M=n+1}^{\infty }%
\frac{\beta _{M+1}}{M}
\end{equation*}%
\begin{equation}
\ll \sum\limits_{k=n+1}^{\infty }\frac{\left\vert a_{k}\right\vert +\beta
_{k}}{k}\ll \sum\limits_{k=n+1}^{\infty }\frac{\beta _{k}}{k}.  \label{s4}
\end{equation}

Let 
\begin{equation*}
\frac{2k\pi }{r}+\frac{\pi }{r}\leq x<\frac{2\left( k+1\right) \pi }{r}-%
\frac{\pi }{N+1},
\end{equation*}%
where $N:=N\left( x\right) \geq r$ and $k=0,1,...,\left[ r/2\right] -1$. By (%
\ref{p3p1}), for $n\geq N,$%
\begin{equation*}
\int\limits_{\frac{2k\pi }{r}+\frac{\pi }{r}}^{\frac{2\left( k+1\right) \pi 
}{r}-\frac{\pi }{n+1}}\left\vert f\left( x\right) -S_{n}\left( f;x\right)
\right\vert dx=\sum\limits_{N=r}^{n}\int\limits_{\frac{2\left( k+1\right)
\pi }{r}-\frac{\pi }{N}}^{\frac{2\left( k+1\right) \pi }{r}-\frac{\pi }{N+1}%
}\left\vert f\left( x\right) -S_{n}\left( f;x\right) \right\vert dx
\end{equation*}%
\begin{equation*}
\ll \sum\limits_{N=r}^{n}\int\limits_{\frac{2\left( k+1\right) \pi }{r}-%
\frac{\pi }{N}}^{\frac{2\left( k+1\right) \pi }{r}-\frac{\pi }{N+1}}\frac{%
\beta _{n+1}}{\left\vert \sin \left( rx/2\right) \right\vert }dx.
\end{equation*}%
Using the inequality%
\begin{equation}
2\left( k+1\right) -\frac{r}{\pi }x\leq \left\vert \sin \left( \frac{r}{2}%
x\right) \right\vert \text{ \ \ for \ \ }x\in \left[ \frac{2k\pi }{r}+\frac{%
\pi }{r},\frac{2\left( k+1\right) \pi }{r}\right]  \label{p3p4}
\end{equation}%
we get%
\begin{equation*}
\int\limits_{\frac{2k\pi }{r}+\frac{\pi }{r}}^{\frac{2\left( k+1\right) \pi 
}{r}-\frac{\pi }{n+1}}\left\vert f\left( x\right) -S_{n}\left( f;x\right)
\right\vert dx\ll \sum\limits_{N=r}^{n}\beta _{n+1}\int\limits_{\frac{%
2\left( k+1\right) \pi }{r}-\frac{\pi }{N}}^{\frac{2\left( k+1\right) \pi }{r%
}-\frac{\pi }{N+1}}\frac{1}{2\left( k+1\right) -\frac{r}{\pi }x}dx
\end{equation*}%
\begin{equation}
=\frac{\pi }{r}\beta _{n+1}\sum\limits_{N=r}^{n}\frac{1}{N}\ll \beta
_{n+1}\ln \left( n+1\right) .  \label{s5}
\end{equation}%
If $N\geq n+1$ then, by (\ref{p3p3}) and (\ref{p3p4}),%
\begin{equation*}
\int\limits_{\frac{2\left( k+1\right) \pi }{r}-\frac{\pi }{n+1}}^{\frac{%
2\left( k+1\right) \pi }{r}}\left\vert f\left( x\right) -S_{n}\left(
f;x\right) \right\vert dx=\sum\limits_{N=n+1}^{\infty }\int\limits_{\frac{%
2\left( k+1\right) \pi }{r}-\frac{\pi }{N}}^{\frac{2\left( k+1\right) \pi }{r%
}-\frac{\pi }{N+1}}\left\vert f\left( x\right) -S_{n}\left( f;x\right)
\right\vert dx
\end{equation*}%
\begin{equation*}
\ll \sum\limits_{N=n+1}^{\infty }\int\limits_{\frac{2\left( k+1\right) \pi 
}{r}-\frac{\pi }{N}}^{\frac{2\left( k+1\right) \pi }{r}-\frac{\pi }{N+1}%
}\left( \sum\limits_{k=n+1}^{N}\left\vert a_{k}\right\vert +\frac{\beta
_{N+1}}{\left\vert \sin \left( rx/2\right) \right\vert }\right) dx
\end{equation*}%
\begin{equation*}
\leq \sum\limits_{N=n+1}^{\infty }\frac{\pi }{N^{2}}\sum%
\limits_{k=n+1}^{N}\left\vert a_{k}\right\vert +\sum\limits_{N=n+1}^{\infty
}\beta _{N+1}\int\limits_{\frac{2\left( k+1\right) \pi }{r}-\frac{\pi }{N}%
}^{\frac{2\left( k+1\right) \pi }{r}-\frac{\pi }{N+1}}\frac{1}{2\left(
k+1\right) -\frac{r}{\pi }x}dx
\end{equation*}%
\begin{equation}
\ll \sum\limits_{k=n+1}^{\infty }\left\vert a_{k}\right\vert
\sum\limits_{N=k}^{\infty }\frac{1}{N^{2}}+\sum\limits_{N=n+1}^{\infty }%
\frac{\beta _{N+1}}{N}\ll \ll \sum\limits_{k=n+1}^{\infty }\frac{\left\vert
a_{k}\right\vert +\beta _{k}}{k}\ll \sum\limits_{k=n+1}^{\infty }\frac{%
\beta _{k}}{k}.  \label{s6}
\end{equation}%
Summing up (\ref{s1}), (\ref{s2}), (\ref{s3}), (\ref{s4}), (\ref{s5}) and (%
\ref{s6}), we finally have%
\begin{equation*}
\left\Vert S_{n}\left( f\right) -f\right\Vert \ll \beta _{n+1}\ln \left(
n+1\right) +\sum\limits_{k=n+1}^{\infty }\frac{\beta _{k}}{k}.
\end{equation*}%
Furthermore by the assumption (\ref{m3}) we obtain that (\ref{m4}) follows
and thus the proof is complete. $\square $

\subsection{Proof of Theorem 5}

If $\left( b_{n}\right) \in RBV\left( \beta ,r\right) $ then%
\begin{equation*}
\sum\limits_{k=1}^{\infty }\frac{\left\vert b_{k}\right\vert }{k}\ll
\sum\limits_{k=1}^{\infty }\frac{\beta _{k}}{k}<\infty .
\end{equation*}%
Thus $b_{n}=o\left( 1\right) $ and by (\ref{l02}) we obtain that for $x\neq 
\frac{2l\pi }{r}$, $l\in 
\mathbb{Z}
$%
\begin{equation*}
\sum\limits_{k=n}^{\infty }a_{k}\cos kx=\frac{1}{2\sin \left( rx/2\right) }%
\left\{ \sum\limits_{k=n}^{\infty }\left( a_{k}-a_{k+r}\right) \cos \left(
k+\frac{r}{2}\right) x\right.
\end{equation*}%
\begin{equation*}
\left. -\sum\limits_{k=n}^{n+r-1}a_{k}\cos \left( k-\frac{r}{2}\right)
x\right\} ,
\end{equation*}%
whence%
\begin{equation*}
\left\vert f\left( x\right) -S_{n}\left( f;x\right) \right\vert =\left\vert
\sum\limits_{k=n+1}^{\infty }a_{k}\cos kx\right\vert
\end{equation*}%
\begin{equation}
\leq \frac{1}{2\left\vert \sin \left( rx/2\right) \right\vert }\left\{
\sum\limits_{k=n+1}^{\infty }\left\vert a_{k}-a_{k+r}\right\vert
+\sum\limits_{k=n+1}^{n+r}\left\vert a_{k}\right\vert \right\} \ll \frac{%
\beta _{n+1}}{\left\vert \sin \left( rx/2\right) \right\vert }.  \label{p1p0}
\end{equation}%
Similarly, as in the previous proof, using (\ref{p1p0}) we can also show that%
\begin{equation*}
\left\Vert S_{n}\left( g\right) -g\right\Vert \ll \beta _{n+1}\ln \left(
n+1\right) +\sum\limits_{k=n+1}^{\infty }\frac{\beta _{k}}{k}.
\end{equation*}%
Thus by the assumption (\ref{m3}) we obtain that (\ref{m5}) follows and our
proof is complete. $\square $

\subsection{Proof of Remark 5}

(i) Let%
\begin{equation*}
a_{n}=\left\{ 
\begin{array}{c}
0\text{ \ \ for \ \ }n\neq 3l+1\text{, }l\in 
\mathbb{N}
\cup \left\{ 0\right\} , \\ 
\frac{1}{\ln n\cdot \ln \left( \ln n\right) }\text{ \ \ for \ \ }n=3l+1\text{%
, }l\in 
\mathbb{N}
.%
\end{array}%
\right.
\end{equation*}%
First we prove that $\left( a_{n}\right) \in RBVS\left( _{5}\beta ,3\right)
. $ If%
\begin{equation*}
A\left( n\right) =\left\{ k:n\leq k\wedge k=3l+1\wedge l\in 
\mathbb{N}
\right\}
\end{equation*}%
then%
\begin{equation*}
\sum\limits_{k=n}^{2n-1}\left\vert a_{k}-a_{k+3}\right\vert
=\sum\limits_{k\in A\left( n\right) }\left\vert a_{k}-a_{k+3}\right\vert
\end{equation*}%
\begin{equation*}
=\sum\limits_{k\in A\left( n\right) }\frac{\ln \left( k+3\right) \cdot \ln
\left( \ln \left( k+3\right) \right) -\ln k\cdot \ln \left( \ln \left(
k\right) \right) }{\ln k\cdot \ln \left( k+3\right) \cdot \ln \left( \ln
k\right) \ln \left( \ln \left( k+3\right) \right) }.
\end{equation*}%
Applying Lagrange's mean value theorem to the function $f\left( x\right)
=\ln x\cdot \ln \left( \ln x\right) $ on the interval $\left[ k,k+3\right] $
we obtain that there exists $z\in \left( k,k+3\right) $ such that%
\begin{equation}
\ln \left( k+3\right) \cdot \ln \left( \ln \left( k+3\right) \right) -\ln
k\cdot \ln \left( \ln \left( k\right) \right) =3\frac{\ln \left( \ln
z\right) +1}{z}.  \label{p2p0}
\end{equation}%
Hence%
\begin{equation*}
\sum\limits_{k=n}^{\infty }\left\vert a_{k}-a_{k+3}\right\vert \leq
3\sum\limits_{k\in A\left( n\right) }\frac{\ln \left( \ln \left( k+3\right)
\right) +1}{k\ln k\cdot \ln \left( k+3\right) \cdot \ln \left( \ln k\right)
\ln \left( \ln \left( k+3\right) \right) }
\end{equation*}%
\begin{equation}
\leq \frac{6}{\ln \left( \ln n\right) }\sum\limits_{k\in A\left( n\right) }%
\frac{1}{k\left( \ln k\right) ^{2}}\ll \frac{1}{\ln n\cdot \ln \left( \ln
n\right) }\ll \sum\limits_{k=\left[ n/c\right] }^{\left[ cn\right] }\frac{%
\left\vert a_{k}\right\vert }{k}.  \label{p2p1}
\end{equation}%
Thus $RBVS\left( _{5}\beta ,3\right) $ and consequently $\left( a_{n}\right)
\in GM\left( _{5}\beta ,3\right) .$

Using (\ref{p3p1}) with $r=3$ and (\ref{p2p1}) we obtain that the series (%
\ref{1}) is convergent for all $x\neq \frac{2l\pi }{3}$, $l\in 
\mathbb{Z}
$, i.e. 
\begin{equation*}
f\left( x\right) =\sum\limits_{k=1}^{\infty }\frac{1}{\ln \left(
3k+1\right) \cdot \ln \left( \ln \left( 3k+1\right) \right) }\cos \left(
3k+1\right) x
\end{equation*}%
for $x\neq \frac{2l\pi }{r}$, $l\in 
\mathbb{Z}
$.

Let $m=\exp \left\{ \exp \left( \exp n\right) \right\} $. Then%
\begin{equation*}
\left\Vert S_{3\left[ m\right] +2}\left( f\right) -S_{3n-1}\left( f\right)
\right\Vert =\frac{1}{\pi }\int\limits_{0}^{\pi }\left\vert
\sum\limits_{k=3n}^{3\left[ m\right] +2}a_{k}\cos kx\right\vert dx
\end{equation*}%
\begin{equation*}
\geq \frac{1}{\pi }\int\limits_{0}^{\frac{2}{3}\pi }\left\vert
\sum\limits_{k=3n}^{3\left[ m\right] +2}a_{k}\cos kx\right\vert dx\geq 
\frac{1}{\pi }\left\vert \sum\limits_{k=3n}^{3\left[ m\right]
+2}a_{k}\int\limits_{0}^{\frac{2}{3}\pi }\cos kxdx\right\vert
\end{equation*}%
\begin{equation*}
=\frac{1}{\pi }\left\vert \sum\limits_{k=3n}^{3\left[ m\right] +2}a_{k}\sin
\left( \frac{2}{3}\pi k\right) \right\vert =\frac{1}{\pi }\left\vert
\sum\limits_{k=n}^{\left[ m\right] }\sum\limits_{l=0}^{2}\frac{a_{3k+l}}{%
3k+l}\sin \left( \frac{2}{3}\left( 3k+l\right) \pi \right) \right\vert
\end{equation*}%
\begin{equation*}
=\frac{1}{\pi }\left\vert \sum\limits_{k=n}^{\left[ m\right]
}\sum\limits_{l=1}^{2}\frac{a_{3k+l}}{3k+l}\sin \left( \frac{2}{3}l\pi
\right) \right\vert =\frac{\sin \frac{2}{3}\pi }{\pi }\sum\limits_{k=n}^{%
\left[ m\right] }\frac{1}{\left( 3k+1\right) \ln \left( 3k+1\right) \cdot
\ln \left( \ln \left( 3k+1\right) \right) }
\end{equation*}%
\begin{equation*}
\gg n-\ln \left( \ln \left( \ln \left( 3n+1\right) \right) \right) .
\end{equation*}%
Since $n-\ln \left( \ln \left( \ln \left( 3n+1\right) \right) \right)
\rightarrow \infty $ as $n\rightarrow \infty $ thus the sequence $%
S_{n}\left( f;x\right) $ does not satisfy the Cauchy condition in $L$-norm.
Hence the sequence $S_{n}\left( f;x\right) $ can not be convergent in $L$%
-norm$.$

(ii) Suppose that%
\begin{equation*}
b_{n}=\left\{ 
\begin{array}{c}
0\text{ \ \ for \ \ }n\neq 2l+1\text{, }l\in 
\mathbb{N}
\cup \left\{ 0\right\} , \\ 
\frac{1}{\ln n\cdot \ln \left( \ln n\right) }\text{ \ \ for \ \ }n=2l+1\text{%
, }l\in 
\mathbb{N}
.%
\end{array}%
\right.
\end{equation*}%
Similarly, as in (i), we can show that $RBVS\left( _{5}\beta ,2\right) $ and
consequently $\left( b_{n}\right) \in GM\left( _{5}\beta ,2\right) .$ It is
clear that the series (\ref{2}) is convergent for $x=l\pi $, $l\in 
\mathbb{Z}
$. Moreover, since $\left( a_{n}\right) \in RBV\left( _{5}\beta ,2\right) $
thus, by (\ref{p1p0}) with $r=2,$ we get that the series (\ref{2}) is
convergent for all $x\neq l\pi $, $l\in 
\mathbb{Z}
$, i.e.%
\begin{equation*}
g\left( x\right) =\sum\limits_{k=1}^{\infty }\frac{1}{\ln \left(
2k+1\right) \cdot \ln \left( \ln \left( 2k+1\right) \right) }\sin \left(
2k+1\right) x\text{ for all }x\in 
\mathbb{R}
.
\end{equation*}%
Let $m=\exp \left\{ \exp \left( \exp n\right) \right\} .$ Similarly as
before we get that%
\begin{equation*}
\left\Vert S_{2\left[ m\right] +1}\left( f\right) -S_{2n-1}\left( f\right)
\right\Vert =\frac{1}{\pi }\int\limits_{0}^{\pi }\left\vert
\sum\limits_{k=2n}^{2\left[ m\right] +1}b_{k}\sin kx\right\vert dx\geq 
\frac{1}{\pi }\left\vert \sum\limits_{k=2n}^{2\left[ m\right]
+1}b_{k}\int\limits_{0}^{\pi }\sin kxdx\right\vert
\end{equation*}%
\begin{equation*}
=\frac{1}{\pi }\left\vert \sum\limits_{k=2n}^{2\left[ m\right]
+1}b_{k}\left\{ 1-\cos \left( k\pi \right) \right\} \right\vert =\frac{1}{%
\pi }\left\vert \sum\limits_{k=n}^{\left[ m\right] }\sum\limits_{l=0}^{1}%
\frac{b_{2k+l}}{2k+l}\left\{ 1-\cos \left( 2k+l\right) \pi \right\}
\right\vert
\end{equation*}%
\begin{equation*}
=\frac{1}{\pi }\left\vert \sum\limits_{k=n}^{\left[ m\right]
}\sum\limits_{l=0}^{1}\frac{b_{2k+l}}{2k+l}\left\{ 1-\cos l\pi \right\}
\right\vert =\frac{2}{\pi }\sum\limits_{k=n}^{\left[ m\right] }\frac{1}{%
\left( 2k+1\right) \ln \left( 2k+1\right) \cdot \ln \left( \ln \left(
2k+1\right) \right) }
\end{equation*}%
\begin{equation*}
\gg n-\ln \left( \ln \left( \ln \left( 2n+1\right) \right) \right) .
\end{equation*}%
Since $n-\ln \left( \ln \left( \ln \left( 2n+1\right) \right) \right)
\rightarrow \infty $ as $n\rightarrow \infty $ thus the sequence $%
S_{n}\left( g;x\right) $ does not satisfy the Cauchy condition in $L$-norm.
Hence the sequence $S_{n}\left( g;x\right) $ can not be convergent in $L$%
-norm$.$ $\square $

\subsection{Proof of Remark 6}

Let $r\geq 3$ and%
\begin{equation*}
d_{n}=\left\{ 
\begin{array}{c}
0\text{ \ \ if \ \ }r\nmid n, \\ 
\frac{1}{n^{2}}\text{ \ \ if \ \ }r\mid n.%
\end{array}%
\right.
\end{equation*}%
It is clear that $\sum\limits_{k=1}^{\infty }\frac{\beta _{k}}{k}<\infty $.
First, we show that $\left( d_{n}\right) \in RBV\left( _{6}\beta ,r\right) $%
. Suppose that $A_{r}\left( n,m\right) =\left\{ k:n\leq k<m\wedge r\mid
k\right\} $. Then%
\begin{equation*}
\sum\limits_{k=n}^{\infty }\left\vert d_{k}-d_{k+r}\right\vert
=\sum\limits_{k\in A_{r}\left( n,\infty \right) }\left\vert
d_{k}-d_{k+r}\right\vert =\sum\limits_{k\in A_{r}\left( n,\infty \right)
}\left\vert \frac{1}{k^{2}}-\frac{1}{\left( k+r\right) ^{2}}\right\vert
\end{equation*}%
\begin{equation*}
=\sum\limits_{k\in A_{r}\left( n,\infty \right) }\frac{r\left( 2k+r\right) 
}{k^{2}\left( k+r\right) ^{2}}\leq 2r\sum\limits_{k\in A_{r}\left( n,\infty
\right) }\frac{1}{k^{3}}\ll \frac{1}{k^{2}}\ll \sum\limits_{k=\left[ n/c%
\right] }^{\left[ cn\right] }\frac{\left\vert d_{k}\right\vert }{k}
\end{equation*}%
and $\left( d_{n}\right) \in RBVS\left( _{5}\beta ,r\right) $. Thus $\left(
d_{n}\right) \in RBVS\left( _{6}\beta ,r\right) .$ Finally, we prove that $%
\left( d_{n}\right) \notin GM\left( _{6}\beta ,2\right) $. Indeed,%
\begin{equation*}
\sum\limits_{k=n}^{2n-1}\left\vert d_{k}-d_{k+2}\right\vert \geq
\sum\limits_{k\in A_{r}\left( n,2n\right) }\left\vert
d_{k}-d_{k+2}\right\vert =\sum\limits_{k\in A_{r}\left( n,2n\right) }\frac{1%
}{k^{2}}\geq \frac{1}{4rn}
\end{equation*}%
and since%
\begin{equation*}
\beta _{n}=\frac{1}{\ln n}\underset{m\geq \left[ n/c\right] }{\max }\left( 
\frac{\ln m}{m}\sum\limits_{v=m}^{2m}\left\vert d_{v}\right\vert \right)
\ll \frac{1}{n}\sum\limits_{v=\left[ n/c\right] }^{\infty }\frac{1}{v^{2}}%
\ll \frac{1}{n^{2}},
\end{equation*}%
the inequality%
\begin{equation*}
\sum\limits_{k=n}^{2n-1}\left\vert d_{k}-d_{k+2}\right\vert \ll \beta _{n}
\end{equation*}%
does not hold, that is, $\left( d_{n}\right) $ does not belong to $GM\left(
_{6}\beta ,2\right) .$ $\square $


\begin{thebibliography}{99}
\bibitem{1} B. Aubertin, J. Fournier,{\it  Integrability theorems for
trigonometric series}, Studia Math., 107 (1) (1993), 33-59.

\bibitem{19} A. S. Belov,{\it  On conditions for the convergence in the mean of
trigonometric Fourier series}, Izv. Tul. Gos. Univ. Ser. Mat. Mekh. Inform.,
4 (1) (1998), 40-46 (in Russian).

\bibitem{20} A. S. Belov, {\it On conditions for convergence (boundedness) in the
mean of partial sums of a trigonometric series}, in: Metric Theory of
Functions and Related Problems in Analysis, Izd. AFT, Moscow, 1999, pp. 1-17
(in Russian).

\bibitem{21} A. S. Belov,{\it  Remarks on the convergence (boundedness) in the
mean of partial sums of a trigonometric series}, Mat. Zametki, 71 (6) (2002),
807-817.

\bibitem{22} S. Fridli,{\it  On the $L_{1}$-convergence of Fourier series}, Studia
Math., 125 (2) (1997), 161-174.

\bibitem{23} J. W. Garrett, C. S. Rees, \v{C}. V. Stanojevi\'{c},{\it  On $L^{1}$
convergence of Fourier series with quasi-monotone coefficients}, Proc. Amer.
Math. Soc., 72 (3) (1978), 535-538.

\bibitem{8} J. W. Garrett, \v{C}. V. Stanojevi\'{c}, {\it Necessary and
sufficient condition for $L^{1}$ convergence of trigonometric series}, Proc.
Amer. Math. Soc., 60 (1976), 68-71.

\bibitem{2} R. J. Le and S. P. Zhou,{\it  A new condition for the uniform
convergence of certain trigonometric series}, Acta Math. Hungar., 108 (1-2)
(2005), 161-169.

\bibitem{9} L. Leindler,{\it  On the uniform convergence and boundedness of a
certain class of sine series}, Anal. Math., 27 (2001), 279-285.

\bibitem{3} E. Liflyand, {\it Lebesgue constants of multiple Fourier series},
Online J. Anal. Com., (1) (2006), available online at http://www.ojac.org/.

\bibitem{16} V. B. Stanojevic,{\it  $L^{1}$-convergence of Fourier series with $O$%
-regularly varying quasimonotonic coefficients}, J. Approx. Theory, 60 (2)
(1990), 168-173.

\bibitem{4} S. A. Teljakovski\u{\i},{\it  On the question of convergence in the
metric $L_{1}$ of the Fourier series with rarely changing Fourier
coefficients}, Approximation of Functions, (2003), 240-243 (in Russian).

\bibitem{5} S. A. Teljakovski\u{\i}, G. A. Fomin,{\it  Convergence in the $L$
metric of Fourier series with quasimonotone coefficients}, Tr. Mat. Inst.
Steklova, 134 (1975), 310-313; translation in Proc. Steklov Inst. Math., 134
(1975), 351-355.

\bibitem{10} S. Tikhonov, {\it Trigonometric series with general monotone
coefficients}, J. Math. Anal. Appl., 326(1) (2007), 721-735.

\bibitem{11} S. Tikhonov, {\it On uniform convergence of trigonometric series},
Mat. Zametki, 81(2) (2007), 304-310, translation in Math. Notes, 81(2)
(2007), 268-274.

\bibitem{12} S. Tikhonov, {\it Best approximation and moduli of smoothness:
Computation and equivalence theorems}, J. Approx. Theory, 153 (2008), 19-39.

\bibitem{24} S. Tikhonov, {\it On $L^{1}$-convergence of Fourier series}, J. Math.
Anal. Appl., 347 (2008), 416-427.

\bibitem{6} Z. Tomovski,{\it  Convergence and integrability for some classes of
trigonometric series}, Dissertationes Math., 420 (2003), 1-65.

\bibitem{13} B. Szal,{\it  A new class of numerical sequences and its
applications to uniform convergence of sine series}, available online at
http://arxiv.org/0905.1294v1.

\bibitem{7} T. F. Xie, S. P. Zhou,{\it  $L^{1}$-approximation of Fourier series
of complex-valued functions}, Proc. Roy. Soc. Edinburgh Sect. A, 126 (2)
(1996), 343-353.

\bibitem{14} D. S. Yu and S. P. Zhou, {\it A generalization of monotonicity and
applications}, Acta Math. Hungar., 115 (3) (2007), 247-267.

\bibitem{17} D. S. Yu, S. P. Zhou and P. Zhou, {\it On $L_{1}$-convergence of
Fourier series under $MVBV$ condition}, available online at
http://arxiv.org/0704.1865v1.

\bibitem{18} D. S. Yu, R. J. Le, S. P. Zhou, {\it Remarks on convergence of
trigonometric series with special varying coefficients}, J. Math. Anal.
Appl., 333 (2) (2007), 1128-1137.

\bibitem{15} A. Zygmund, {\it Trigonometric series}, Vol. I, University Press
(Cambridge, 1959).
\end{thebibliography}
\end{document}